\documentclass{amsart}
\usepackage{graphics}

\def\ba{\begin{array}}

\def\ba{\begin{array}}
\def\ea{\end{array}}

\def\d{\delta}
\def\g{\gamma}
\def\l{\lambda}
\def\q{\quad}
\def\ep{\varepsilon}
\def\f{\frac}

\def\l2{L^2(0,1)}

\def\w2p{W^{2,p}(0,1)}

\def\N{\mathcal N}
\def\R{\mathcal R}

\begin{document}
\begin{center}
{\Large\bf On deconvolution problems: numerical aspects.}
\vskip7mm

Alexander G. Ramm\\
E-mail: ramm@math.ksu.edu\\
Department of Mathematics\\
Kansas State University\\Manhattan, KS 66506, USA.\\
\vspace{0.5cm}

Alexandra B. Smirnova\\
E-mail: smirn@mathstat.gsu.edu\\
Department of Mathematics and Statistics\\
Georgia State University\\Atlanta, GA 30303, USA.\\
\end{center}
\vskip 5mm


\vspace{0.5cm}

{\bf Abstract.} An optimal algorithm is described for solving the
deconvolution problem of the form ${\bf k}u:=\int_0^tk(t-s)u(s)ds=f(t)$
given the noisy data $f_\delta$, $||f-f_\delta||\leq \delta.$
The idea of the method consists of the representation ${\bf k}=A(I+S)$,
where $S$ is a compact operator, $I+S$ is injective, $I$ is the identity
operator, $A$ is not boundedly invertible, and an optimal regularizer is
constructed for $A$. The optimal regularizer is constructed using
the results of the paper MR 40\#5130.

\section{{\bf Introduction}}
\setcounter{equation}{0} \setcounter{theorem}{0}
\renewcommand{\thetheorem}{1.\arabic{theorem}}
\renewcommand{\theequation}{1.\arabic{equation}}
\vspace*{-0.5pt}

Deconvolution problem consists of solving equation of the form
\begin{equation}\label{1.1}
{\bf k}u:=\int^t_{0}k(t-s)u(s)\,ds:=k\star u=f(t),\q 0\le t\le T,
\end{equation}
where $k(t)$, $t\ge 0,$ is a kernel of linear integral equation
(\ref{1.1}), $k\star u$ is the convolution. It is important in
many engineering applications, in physics, and other areas. There
is a vast literature on deconvolution methods, see, for example,
\cite{gv}.

If the operator ${\bf k}$ in (\ref{1.1}) is considered as an
operator on $X:=L^\infty(0,T)$, and $\int^T_0|k(t)|\,dt< \infty$,
then ${\bf k}$ is not boundedly invertible, so problem (\ref{1.1})
is ill-posed. Assume that  the data $f$ are noisy:
$f_\d$ is given, such that $||f-f_\d||\le \d$. In
this case it is natural to seek an approximate solution  of
equation (\ref{1.1}) in the class $Q_\d:=\{u\in X: \,\,||{\bf
k}u-f_\d||\le \d \}$. However, for ill-posed equation (\ref{1.1})
an arbitrary element $u_\d\in Q_\d$ cannot be taken as an
approximate solution to (\ref{1.1}), since $u_\d$ is not
continuous with respect to $\d$ in general. In order to select
possible solutions one needs to use {\it a priori} information
(usually available) about the solution, which may be of a
quantitative or qualitative nature.

The usage of qualitative {\it a priori} information makes it
possible to narrow the class of solutions, for example, to a
compact set, so that the problem becomes stable under small
changes in the data. This leads to a concept of a {\bf
quasisolution} \cite{ivt}. Various algorithms for
approximate determination of quasisolutions were studied in
\cite{ivt}.

{\it A priori} information of a qualitative nature (for example,
smoothness of the solution) generates different approaches. The
one which is used often is  {\bf variational regularization}
\cite{ta}, \cite{ph}, which allows one to construct stable
approximate solutions to ill-posed problems by means of a
stabilizing functional. The variational method has been
extensively developed in \cite{gr}, \cite{ehn}, and certain {\it a
priori} and {\it a posteriori} choices of a regularization
parameter $\ep=\ep(\d)$ have been designed and implemented
\cite{m},\cite{en}.

 One can also find approximate solutions to (\ref{1.1}) {\bf by
 iterations}
(see \cite{vv}, \cite{bg}), taking
$x_n=R(f_\d,x_{n-1},...,x_{n-k})$, where $k\le n$. For these
solutions to be stable under small changes of the data,
the iteration number $n=n(\delta)$ yielding $x_n$ must
depend on the $\delta$ suitably.

Other important techniques  in theory of ill-posed problems give
regularizing operators by using Fourier, Laplace, Mellin, and
other integral transforms, statistical regularization, and the
dynamical systems method (DSM) \cite{r451}, \cite{r459}).

In \cite{rg} some general new approaches are proposed for solving
an ill-posed  deconvolution problem. One of these approaches is
based on the following idea. Assume that the operator ${\bf k}$ in
(\ref{1.1}) can be decomposed into a sum ${\bf k}:=A+B$, where the
operator $A^{-1}B:=S$ is compact in the Banach space $X$, in which
${\bf k}$ acts, and $I+S$ is boundedly invertible. By the Fredholm
alternative, it is equivalent to assuming that $\N(I+S)=\{0\}$,
where $\N(A)$ is the null space of $A$. In this case $I+S$ is an
isomorphism of $X$ onto $X$, $\R(A)=\R({\bf k})$, where
$\R(A)$ is the range of the operator $A$, and
\begin{equation}\label{1.2}
{\bf k}u=A(I+S)u=f_\d.
\end{equation}
If a regularizer for $A$ is known, then (\ref{1.2}) can be solved
stably by the scheme
\begin{equation}\label{1.3}
u_\d=(I+S)^{-1}R(\d)f_\d,
\end{equation}
and
\begin{equation}\label{1.4}
||u-u_\d||\to 0 \q \mbox{as} \q \d\to 0.
\end{equation}
Since $I+S$ is an isomorphism, the error $||v-v_\d||$ of the
approximation of the solution of the equation $Av=f_\d$ by the
formula $v_\d=R(\d)f_\d$ is of the same order as $||u_\d-u||$. In
this paper (see sections 2 and 3) we show that the proposed method
is practically efficient and works better than the variational
regularization.

Theoretically the proposed method is optimal on the class of
the data defined as a triple $\{\delta, f_\delta, M_2\}$,
where
$f\in C^2(0,T)$,  $||f^{\prime\prime}||\leq
M_2$, and $f$ is otherwise arbitrary,
$f_\delta\in L^\infty(0,T)$ and $||f-f_\delta||\leq
\delta$ and $f_\delta$ is otherwise arbitrary.

The operator $R(\delta)$, defined in (2.3) and originally
proposed in [10] for
stable numerical differentiation, yields
an optimal estimate of $f'$ in $L^\infty(0,T)-$norm in the
following sense:
$$\inf_T \sup_{\{f_\delta:\, ||f-f_\delta||\leq \delta,\,\,
||f||\leq M_2\}} ||Tf_\delta-f'||\geq (2M_2\delta)^{1/2},
$$
where the $infimum$ is taken over all, linear and
non-linear, operators $T:X\to X$, $X=L^\infty(0,T)$, the
$supremum$ is taken over all $f$ and $f_\delta$ subject to
the conditions $f\in C^2(0,T)$,  $||f^{\prime\prime}||\leq
M_2$, $||f-f_\delta||\leq \delta$,
and
$$
||R(\delta)f_\delta-f'||\leq (2M_2\delta)^{1/2},
$$
(see e.g., [16],[17], [14]).

This argument shows that our "deconvolution" method
for stable solution of (1.1) is optimal on the above data
set: the operator $R(\delta)$ gives an optimal (on the above
data set) approximation of $f'$. Inversion of an
isomorphism $I+S$, where $S$ is a compact operator, can be
done very accurately by a projection method, for example, so
that the total error of the solution is of the same order
as the error obtained by applying $R(\delta)$.

\vskip 5mm

\section{\bf The case $k(t)\in C^1(0,T)$}
\setcounter{equation}{0} \setcounter{theorem}{0}
\renewcommand{\thetheorem}{2.\arabic{theorem}}
\renewcommand{\theequation}{2.\arabic{equation}}
\vskip5mm

\vspace{0.5cm}

\begin{figure}[t]
\begin{center}
   \scalebox{1.3}{\resizebox{90mm}{80mm}{{
      \includegraphics{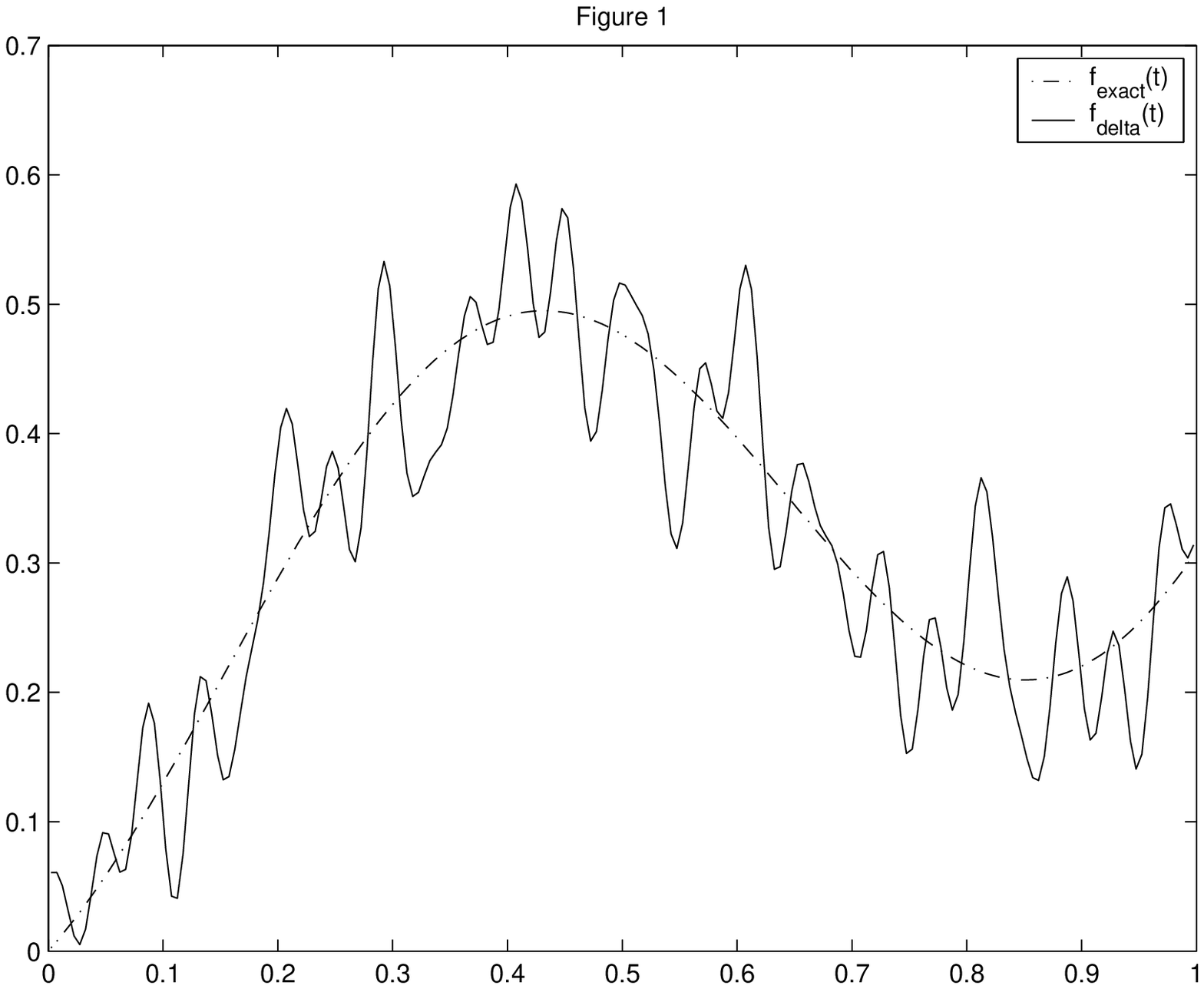}}}}
\end{center}
\end{figure}

Let  $k(t)\in C^1(0,T)$ and $k(0)\neq 0$. Then without loss of
generality one can take $k(0)=1$. As in \cite{rg}, write
(\ref{1.1}) as
\begin{equation}\label{2.2}
{\bf k}u=\int^t_{0}u(s)\,ds+\int^t_{0}[k(t-s)-1]u(s):=Au+Bu=f.
\end{equation}
Assume that $f(x)$ is given by its $\d$-approximation, i.e.
one knows $f_{\d}(x)$ such that $||f-f_\d||_X\le
\d$. In the experiments of this section $\d=0.1$. Let
$A^{-1}B:=S$. Then
\begin{equation}\label{2.3}
{\bf k}u=A(I+S)u=f.
\end{equation}
Stable inversion of $A$ is equivalent to stable numerical
differentiation of noisy data, and therefore as a regularizer
$R(\d)f_\d$ for $A$ one can use (see
\cite{r1},\cite{r397},\cite{r133}, \cite{r71}, and also
\cite{r415},\cite{r441})
\begin{equation}\label{2.4}
R(\d)f_\d:=\f{f_\d(t+h(\d))-f_\d(t-h(\d))}{2h(\d)},
\end{equation}
with $h(\d)=\left(\f{2\d}{M_2}\right)^{1/2},\q
||f''||_{L^\infty_{(0,T)}}\le M_2$. Hence
\begin{equation}\label{2.5}
(I+S)u_\d=R(\d)f_\d,
\end{equation}
where $S$ is a Volterra operator:
$Su_\d=\int^t_{0}k'(t-s)u_\d(s)\,ds.$ To test numerical
efficiency of the above deconvolution algorithm, we take
\begin{equation}\label{2.6}
k(y)=\exp(ay),\q
f(t)=\f{(b+a)(\exp(at)-\cos(bt))+(b-a)\sin(bt)}{a^2+b^2}.
\end{equation}
Then equation (\ref{1.1}) has the
exact solution:
\begin{equation}\label{2.7}
u_{orig}(t)=\sin(bt)+\cos(bt).
\end{equation}
The graphs of  $f$ and its $\d$-approximation, $f_\d$, for
$T=1,\,$ $\, a=1,\,$ $\, b=2\pi,\,$ are presented in Figure 1. The
perturbation was generated as a sum of five sinusoids with various
periods and amplitudes in such a way that $\,||f-f_{\d}||_X\le
0.1.\,$ For $\d=0.1$ and for the above choice of $f$, $T$, $a$,
and $b$, one has $h(\d)=\left(\f{2\d}{M_2}\right)^{1/2}=0.1253$.
Since in practice often only an estimate for $M_2$ may be
available, our first experiment was done with the approximate
value of $h(\d)$, namely $h=0.105$. The goal of the first
experiment was to compare the results obtained by the
deconvolution method suggested in \cite{rg} and by the variational
regularization with a choice of the parameter by the Morozov
discrepancy principle. The integral in (\ref{1.1}) was calculated
by the corrected trapezoid formula (see \cite{dr})
with
the number of node points $n=200$ on the interval $\,[0,1].\,$ The
graphs of $u_{\mbox{disc}}(t)$ and $u_{\mbox{deconv}}(t)$ as well
as the graph of the original solution, $u_{\mbox{orig}}(t)$, for
$h(\d) = 0.105$ and $n=200$ are given in Figure 2. One can see
from the picture that method \cite{rg} provides higher quality of
reconstruction.

\vspace{0.5cm}

\begin{figure}[t]
\begin{center}
   \scalebox{1.3}{\resizebox{90mm}{80mm}{{
      \includegraphics{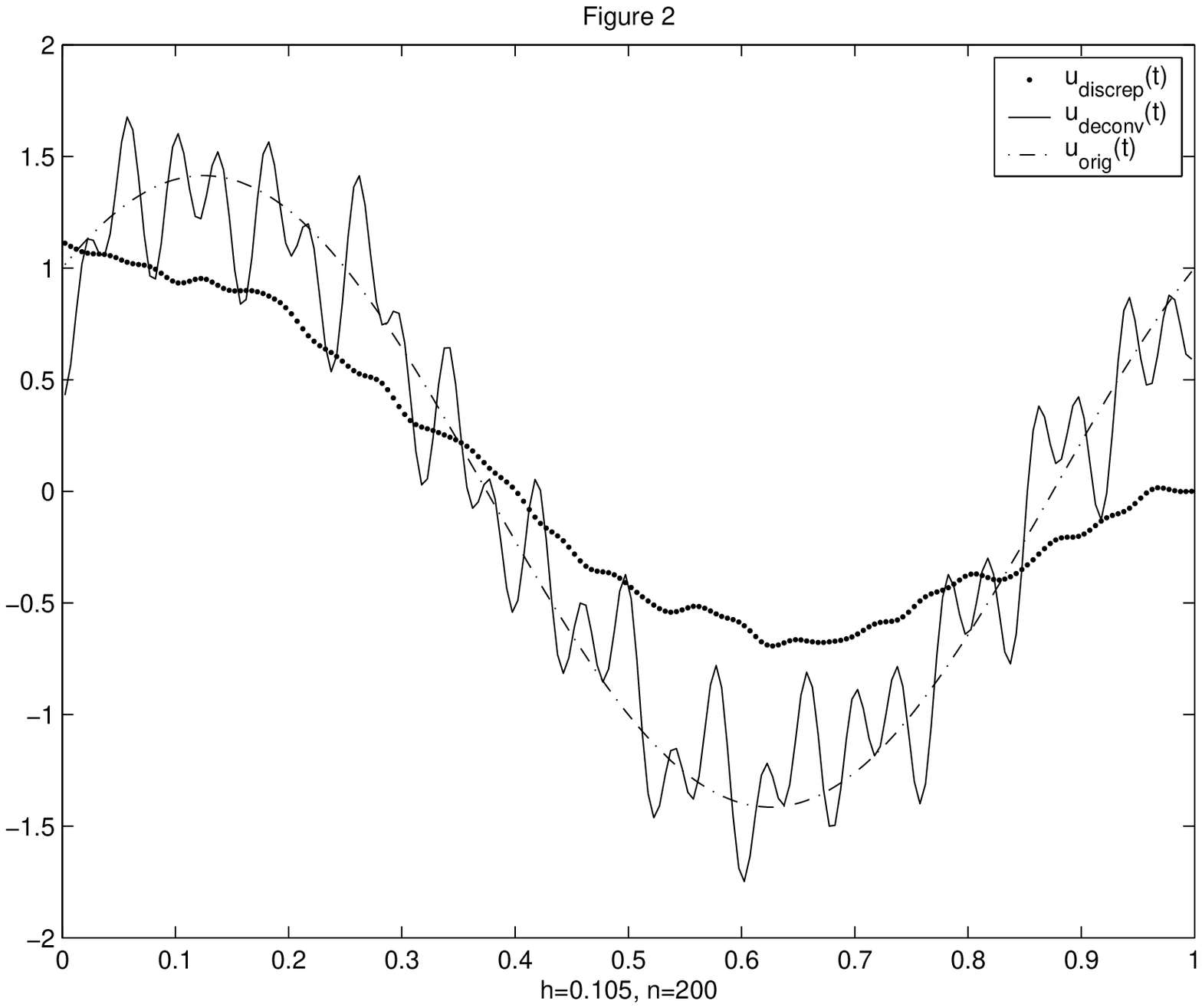}}}}
\end{center}
\end{figure}

\vspace{0.5cm}

  \centerline{Table 1.}

\vspace{0.5cm}

\begin{tabular}{|c|c|c|c|}
\hline
&&&\\
$t$ & $\hskip 1cm u_{\mbox{exact}}(t)\hskip 1cm$ & $\hskip 1cm u_{\mbox{disc}}(t)\hskip 1cm$
& $\hskip 1cm u_{\mbox{deconv}}(t)\hskip 1cm$\\
&&&\\
\hline
&&&\\
0.05 &1.26007351067010&   0.88613219081253&   1.61104047434242\\
0.15&1.39680224666742&   0.77345683250358&   1.16771020714854\\
0.25&   1.00000000000000&   0.78531546607804 &  0.97567292365993\\
0.35&   0.22123174208247&   0.32264819143761&   0.46901890046136\\
0.45&  -0.64203952192021&  -0.01522580641369 & -0.94010917284100\\
0.55&  -1.26007351067010 & -0.72058597578420 & -1.39931254313538\\
0.65&  -1.39680224666742 & -0.65363525334725 & -1.13246945274454\\
0.75&  -1.00000000000000 & -0.84181827797783 & -1.26012127085008\\
0.85&  -0.22123174208247 & -0.48659287989254 & -0.24854842261471\\
0.95&   0.64203952192021 & -0.25478764331776 &  0.99713489435843\\
\hline
\end{tabular}

\vspace{0.7cm}

Table 1 allows one to analyze the computed values of
$u_{\mbox{disc}}(t)$ and $u_{\mbox{deconv}}(t)$ for $h(\d)=0.1$
and $n=10$. The regularization parameter for the variational
regularization calculated by the Morozov discrepancy principle,
$\ep_{\mbox{disc}}$, is equal to $0.0275$ for our particular
$f_\d$. The functions $\,u_{\mbox{disc}}(t)$ and
$\,u_{\mbox{deconv}}(t)$ approximate the exact solution
$u_{\mbox{orig}}(t)$ with the relative errors ${\bf
\delta_{\mbox{disc}}=0.5216}$ and ${\bf
\delta_{\mbox{deconv}}=0.2470}$, respectively, for $n=10$.

\vspace{0.5cm}

\begin{figure}[t]
\begin{center}
   \scalebox{1.3}{\resizebox{90mm}{80mm}{{
      \includegraphics{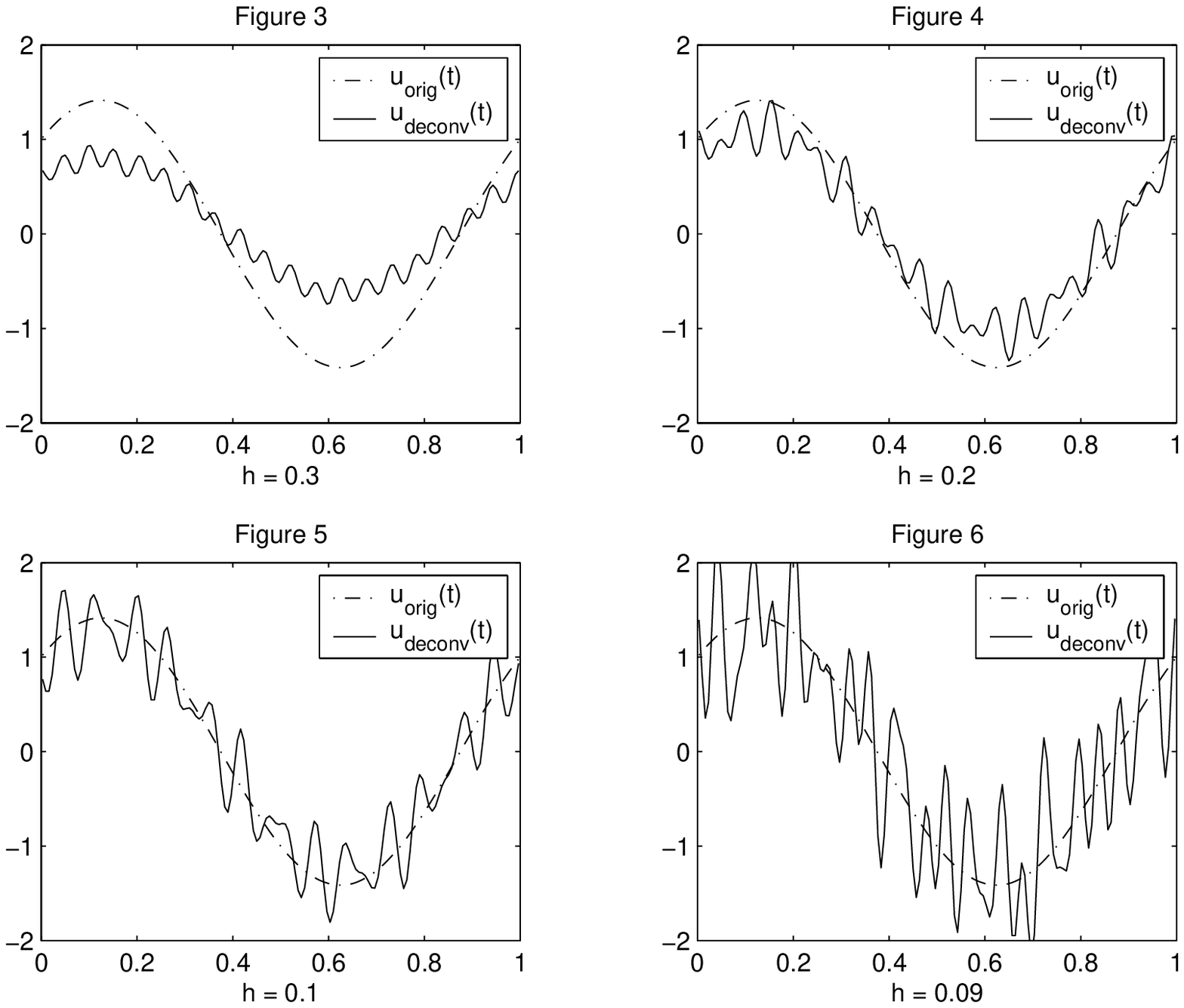}}}}
\end{center}
\end{figure}


  \centerline{Table 2.}

\vspace{0.5cm}

\begin{tabular}{|c|c|c|}
\hline
&&\\
$n$& $\delta_{\mbox{disc}}$& $\delta_{\mbox{deconv}}$ \\
&&\\
\hline
&&\\
10&0.52160739359373&0.24703402714545\\
&&\\

50&0.47882066139400& 0.29886579484582\\
&&\\
100& 0.49421933901812&0.29887532874922
\\
&&\\
\hline
\end{tabular}

\vspace{0.5cm}

The deconvolution procedure \cite{rg}
is  applicable when the constant $M_a$, $a>1$ is known.
 Here $M_a$ is the bound on the $f^{(a)}$, $a>0$ is a real
number, and $f^{(a)}$ is
the (fractional order) derivative of $f$ (see \cite{r441}
for details).
Figures 3-6 show the dependence of the quality of calculations
provided by the deconvolution technique for different values of
$h(\d)$ with the same $f_\d$ that is given in Figure 1. The level
of reconstruction is acceptable for all values of $h(\d)\in
(0.09,0.3)$, but the best quality is attained for the near-optimal
values: $h(\d)=0.1$ and $h(\d)=0.2$. Outside the interval $(0.09,
0.3)$ the reconstruction by the variational regularization works
better because $h(\d)$ is away from its optimal value.

Table 2 contains relative errors, $\delta_{\mbox{disc}}$ and
$\delta_{\mbox{deconv}}$, for values of $n=10,\,\,50,\,\,100$.  In
both cases the relative errors are  not decaying further as $n$
increases, because the major component in these errors come from
the noise level, and not from the error of the computational
methods.

\vskip 7mm
\section{\bf Kernel of the type
$\,k(t)=\f{t^{\g-1}}{\Gamma(\g)}+m(t),\q 0<\g<1,\q m(t)\in C^1$}
\setcounter{equation}{0} \setcounter{theorem}{0}
\renewcommand{\thetheorem}{3.\arabic{theorem}}
\renewcommand{\theequation}{3.\arabic{equation}}

\begin{figure}[t]
\begin{center}
   \scalebox{1.3}{\resizebox{90mm}{80mm}{{
      \includegraphics{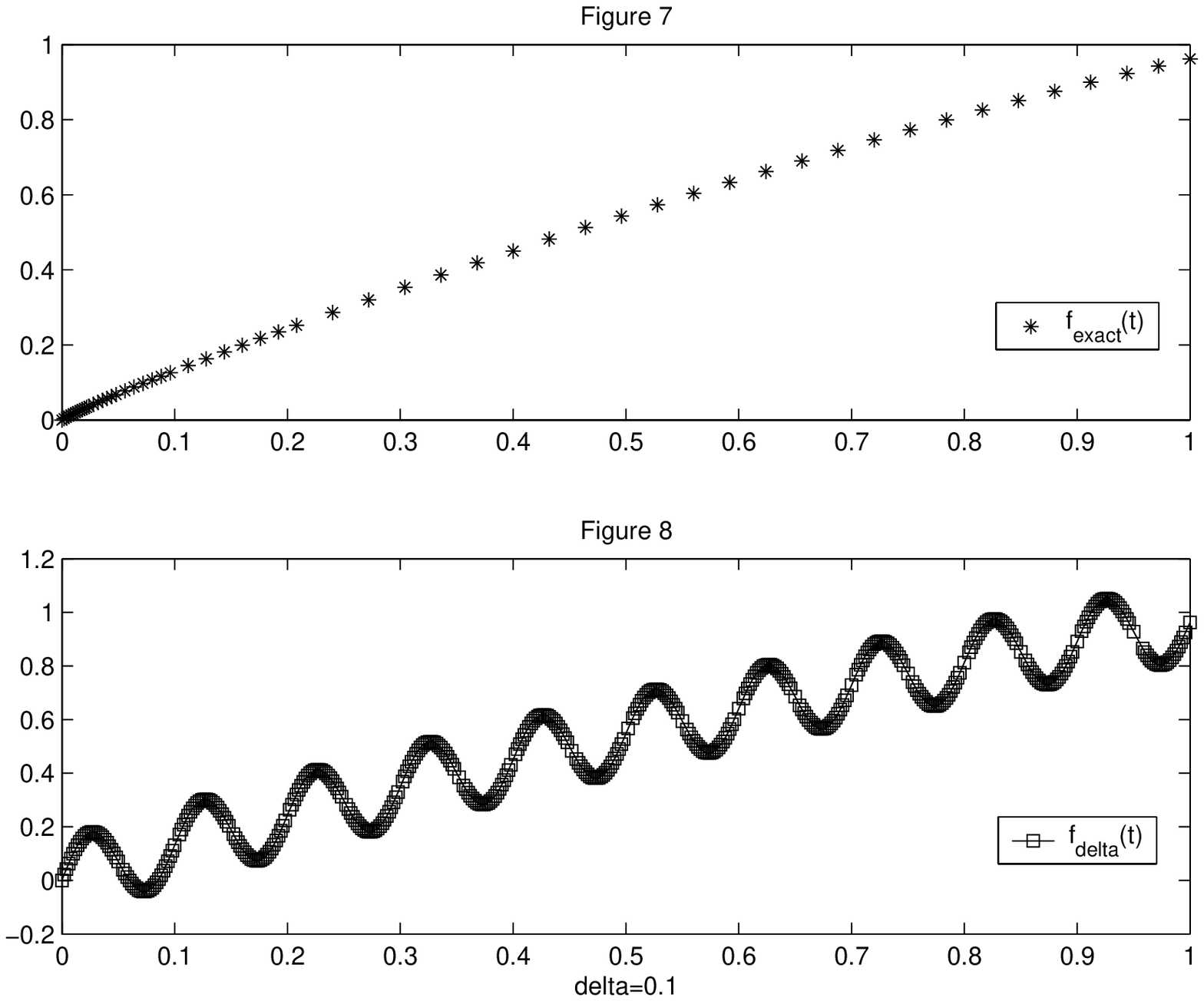}}}}
\end{center}
\end{figure}

\vspace{0.7cm}

In this section we solve  (\ref{1.1}) with the kernel $k(t)$  of
the form
\begin{equation}\label{3.1}
k(t)=\f{t^{\g-1}}{\Gamma(\g)}+m(t),\q0<\g<1,\q m(t)\in C^1
\end{equation}
As in  \cite{rg}, write (\ref{1.1}) as
\begin{equation}\label{3.2}
{\bf k}u:=Au+Bu=f,
\end{equation}
where
\begin{equation}\label{3.3}
Au:=\f{t^{\g-1}}{\Gamma(\g)}\star u,\q Bu:=m\star u.
\end{equation}
One has (\cite{gsh}, pp.117-118)
$\,A^{-1}f=\f{1}{\Gamma(1-\g)}\int^t_0\f{f'(s)}{(t-s)^{\g}}\,ds.\,$
Since the right-hand side $f$ is given by its $\d$-approximation
$f_\d,$ $\,||f-f_\d||_X\le \d$, we replace $A^{-1}$ by the
regularizer $R_1(\d)$ (see \cite{rg}):
\begin{equation}\label{3.4}
R_1(\d)f_\d:=\f{1}{\Gamma(1-\g)}\int^t_0\f{(R(\d)f_\d)(s)}{(t-s)^{\g}}\,ds.
\end{equation}
The operator $R(\d)$ in (\ref{3.4}) is defined by formula
(\ref{2.4}) with $h={\bf 0.12}$. One gets
\begin{equation}\label{3.5}
(I+S)u_\d=R_1(\d)f_\d.
\end{equation}
and
\begin{equation}\label{3.6}
 Su_\d:=A^{-1}Bu_\d=\f{1}{\Gamma(1-\g)}\int^t_0\f{m(0)u_\d(s)+\int^s_0
 m'(s-p)u_\d(p)dp}{(t-s)^{\g}}\,ds.
\end{equation}
The goal of the experiment was to compare two numerical methods
for solving (\ref{1.1})-(\ref{3.1}): deconvolution method
(\ref{3.5})-(\ref{3.6}) and variational regularization with a
choice of the parameter by the discrepancy principle.

The function
$$f(t)=\f{t^{\g}}{\Gamma(1+\g)}\left(1-\f{2t^2}{(1+\g)(2+\g)}\right)
+\f{t^3}{3}\left(1-\f{t^2}{10}\right),\q t\in [0,1],$$ was chosen
as the solution to direct problem (\ref{1.1})-(\ref{3.1}) with
$m(t)=t^2$ and the model function $u_{exact}(t)=1-t^2$. Then for
the numerical tests the noisy function $f_\d$, $\,||f-f_\d||_X\le
\d$, $\d=0.1,$ was used. The graphs of $f$ and $f_\d$ for $\g=0.1$
are given in Figures 7 and 8.

\vspace{0.5cm}

\begin{figure}[t]
\begin{center}
   \scalebox{1.3}{\resizebox{90mm}{80mm}{{
      \includegraphics{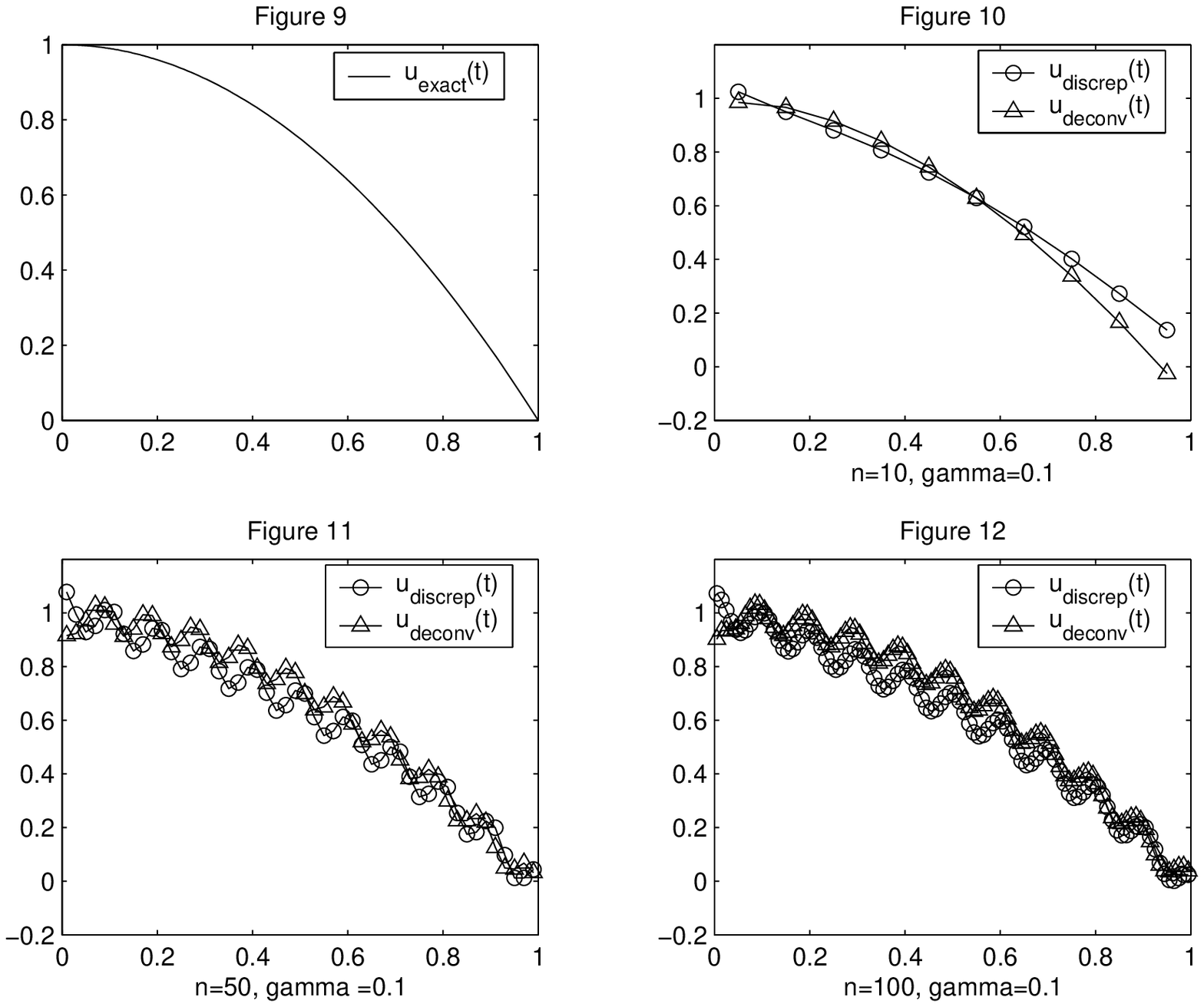}}}}
\end{center}
\end{figure}

\vspace{0.5cm}

  \centerline{Table 3.}

\vspace{0.5cm}

\begin{tabular}{|c|c|c|c|}
\hline
&&&\\
$t$ & $\hskip 1cm u_{\mbox{exact}}(t)\hskip 1cm$ & $\hskip 1cm
u_{\mbox{disc}}(t)\hskip 1cm$
& $\hskip 1cm u_{\mbox{deconv}}(t)\hskip 1cm$\\
&&&\\
\hline
&&&\\
0.05&   0.99000000000000  & 0.93361656127658 &  1.00281244943820\\
0.15&   0.96000000000000 &  0.85983757148008 &  0.98540317766041\\
0.25&   0.91000000000000 &  0.78772680932067 &  0.93442776030698\\
0.35&   0.84000000000000 &  0.71362820985483 &  0.85896131974899\\
0.45&   0.75000000000000 &  0.63504318796224 &  0.76170936024357\\
0.55&   0.64000000000000 &  0.55028612377340 &  0.64396777696250\\
0.65&   0.51000000000000 &  0.45824705747747 &  0.50654835394340\\
0.75&   0.36000000000000 &  0.35823345537370 &  0.35004574315540\\
0.85&   0.19000000000000 &  0.24979168725846 &  0.17493299888351\\
0.95&                  0 &  0.13214536398250  &-0.01840021170081\\
\hline
\end{tabular}

\vspace{0.7cm}


Figures 10-12 illustrate the numerical performance of method
(\ref{3.5})-(\ref{3.6}) and variational regularization for
$\g=0.1$. The solutions evaluated by formulas
(\ref{3.5})-(\ref{3.6}) and by variational regularization for
$n=10$ and $\g=0.1$ are also presented in Table 3.

\begin{figure}[t]
\begin{center}
   \scalebox{1.3}{\resizebox{90mm}{80mm}{{
      \includegraphics{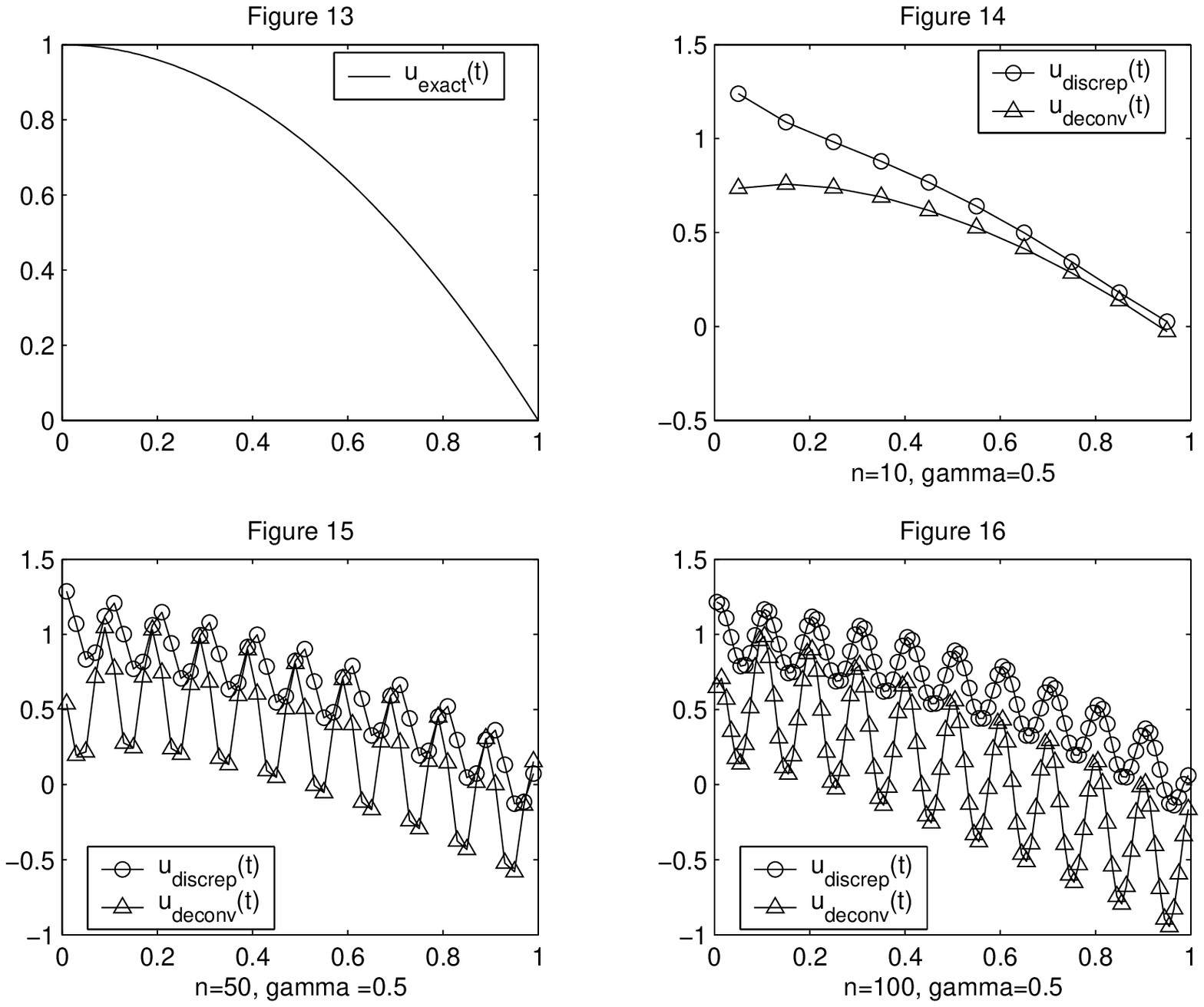}}}}
\end{center}
\end{figure}

\vspace{0.5cm}


\begin{figure}[t]
\begin{center}
   \scalebox{1.3}{\resizebox{90mm}{80mm}{{
      \includegraphics{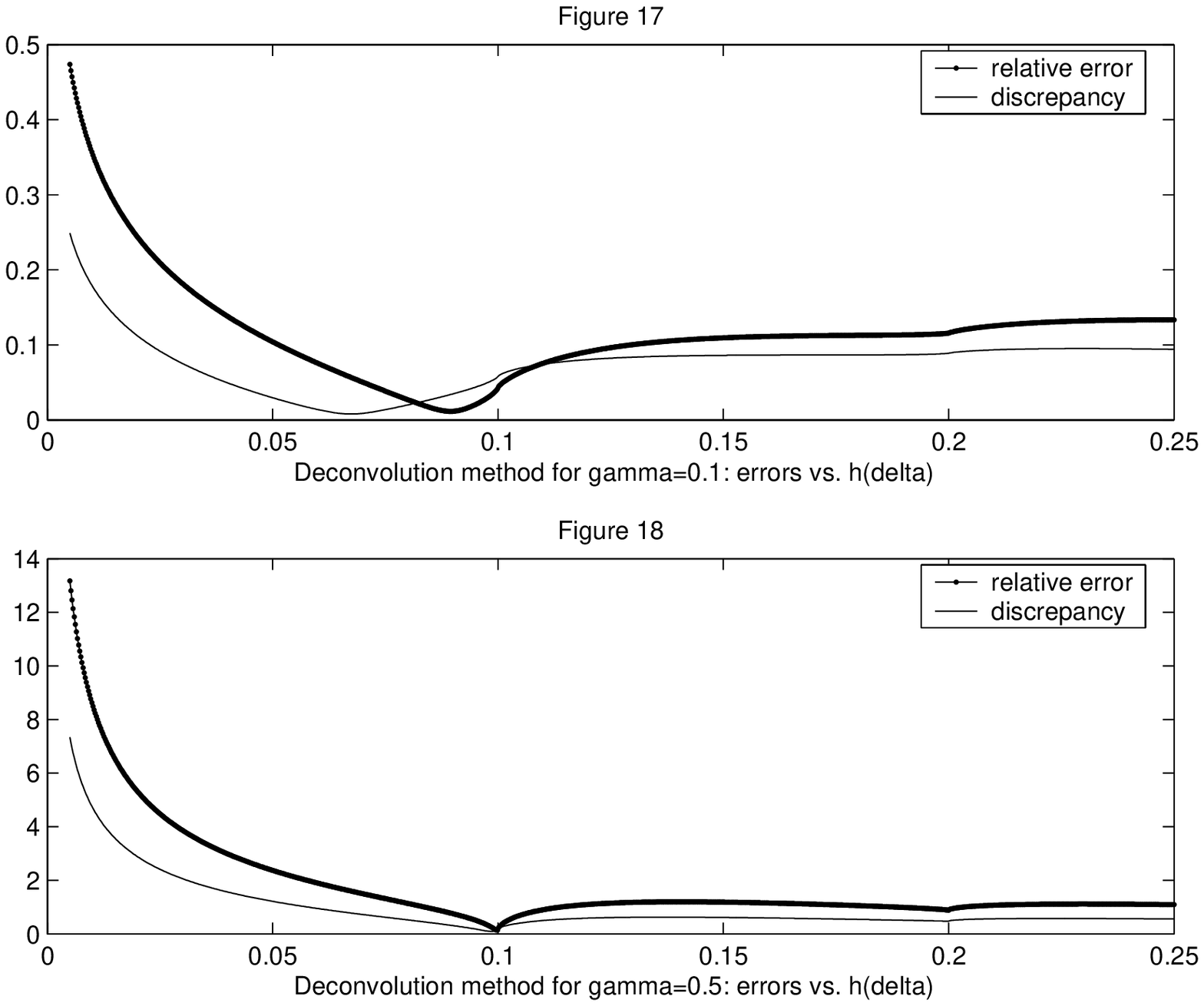}}}}
\end{center}
\end{figure}







The results obtained for our particular test problem show that for
small values of $\g$ the deconvolution approach is superior to
variational regularization both in terms of accuracy and
stability. However as $\g$ is getting bigger, the efficiency of
the deconvolution method (as well as the efficiency of variational
regularization) is getting worse. This is happening because when
$\gamma$ is close to $1$, the ill-posedness of problem
(\ref{1.1})-(\ref{3.1}) grows due to the errors in calculations of
the singular integral. One can compare Figures 9-12 and 13-16.
Moreover, as $\g$ changes from $0.1$ to $0.9$, method
(\ref{3.5})-(\ref{3.6}) becomes very sensitive  to slight
variations of $h(\d)$. To illustrate this phenomena,  we present
the dependence of relative errors and discrepancies on $h(\d)$ for
$\g=0.1$ and $\g=0.5$ in Figures 17 and 18. For $\g=0.1$ the
relative error of the deconvolution method remains less than
$10\%$ when $h(\d)\in (0.05, 0.12)$, while for $g=0.5$ the
relative error is only small for $h=0.1$.

Finally, it is important to mention that CPU time for both
methods, (\ref{3.5})-(\ref{3.6}) and variational regularization,
is approximately the same and it is very small: about 3-4
milliseconds for $n=200$.

{\bf Conclusion.} The paper presents numerical results of the
implementation of the deconvolution method developed by AGR and
presented together with other results  in [14]. The method
is shown to be optimal in the sense explained in Section 1.
The numerical results confirm the theoretical results on which the method
is based. It is shown that the method is more accurate than the
variational regularization method with the regularization parameter chosen
by the discrepancy method.

\vspace{0.5cm}

\bibliographystyle{amsplain}

\end{document}